\documentclass[12pt,psamsfonts]{amsart}

\newtheorem{anyprop}{Anyprop}[section]

\newtheorem{theorem}[anyprop]{Theorem}

\newtheorem{proposition}[anyprop]{Proposition}

\theoremstyle{remark}

\newtheorem{problem}[anyprop]{Problem}

\newtheorem{remark}[anyprop]{Remark}

\newcommand{\NN}{\mathbb{N}}
\newcommand{\ZZ}{\mathbb{Z}}

\newcommand{\FF}{\mathbb{F}}

\newcommand{\PP}{\mathbb{P}}

\newcommand  {\shF}     {\mathcal{F}}

\newcommand  {\shL}     {\mathcal{L}}

\newcommand  {\shQ}     {\mathcal{Q}}

\newcommand  {\shS}     {\mathcal{S}}
\newcommand  {\shT}     {\mathcal{T}}

\newcommand  {\foa}     {\mathfrak{a}}

\newcommand  {\fom}     {\mathfrak{m}}

\newcommand  {\dual}    {\vee}

\newcommand  {\Ext}     {\operatorname{Ext}}

\newcommand  {\GL}      {\operatorname{GL}}

\newcommand  {\lra}     {\longrightarrow}

\renewcommand{\O}       {\mathcal{O}}

\newcommand  {\Proj}    {\operatorname{Proj}}

\newcommand  {\ra}      {\rightarrow}

\newcommand  {\rk}    {\operatorname{rk}}

\newcommand  {\Spec}    {\operatorname{Spec}}

\newcommand  {\Syz}     {\operatorname{Syz}}

\newcommand {\gr} {+ \rm gr}

\newcommand{\stacklra}[1]{ \stackrel{ #1 }{\lra} }
\newcommand{\comdots}{ , \ldots , }
\newcommand{\komdots}{ , \ldots , }

\newcommand{\subsetdots}{ \subset \ldots \subset }

\newcommand{\numiii}{\renewcommand{\labelenumi}{(\roman{enumi})}}

\numberwithin{equation}{section}

\usepackage{amscd}
\usepackage{amssymb}

\def\mydate{\number\day\space\ifcase\month \or January\or February\or March\or April\or May\or
June\or July\or August\or September\or October\or November\or
December\fi \space\number\year}

\begin{document}

\title[Tight closure and plus closure]
{Tight closure and plus closure in dimension two}

\author[Holger Brenner]{Holger Brenner}
\address{Department of Pure Mathematics, University of Sheffield,
  Hicks Building, Hounsfield Road, Sheffield S3 7RH, United Kingdom}
\email{H.Brenner@sheffield.ac.uk}

\subjclass{}

\begin{abstract}
We prove that the tight closure and the graded plus closure of a homogeneous
ideal coincide for a two-dimensional $\NN$-graded domain
of finite type over the algebraic closure of a finite field. This answers in this case a ``tantalizing question'' of Hochster.
\end{abstract}

\maketitle

\noindent
Mathematical Subject Classification (2000):
13A35; 14H60

\section*{Introduction}

One of the basic open questions in tight closure theory is the
problem whether the tight closure of an ideal $I$ in a Noetherian
domain $R$ of positive characteristic $p$ is the same as its plus
closure. Hochster calls this a tantalizing question \cite[Remark to
Theorem 3.1]{hochstertightsolid}. A positive answer to this question
would also give a positive answer to the localization problem.

The tight closure of an ideal $I =(f_1 , \ldots ,f_n) \subseteq R$
is the ideal given by
$$I^* = \{f \in R:\, \exists c \neq 0 \mbox{ such that }
cf^q \in (f_1^q, \ldots ,f_n^q) \mbox{ for all powers } q=p^{e}\} \,
.$$ The plus closure of $I$ is just $I^+=R \cap IR^+$, where $R^+$
is the integral closure of $R$ in an algebraic closure of the
quotient field $Q(R)$. An element $f$ belongs to this plus closure
if and only if there exists a finite extension $R \subseteq S$ such
that $f \in IS$. The inclusion $I^+ \subseteq I^*$ is easy, see
\cite[Theorem 1.7]{hunekeapplication}.

K. Smith proved in \cite{smithparameter}
that the answer to this question is yes
for para\-meter ideals in an excellent local domain.
Recall that $d$ elements in a local ring $(R, \fom)$ of dimension $d$
are called parameters if they generate
an ideal primary to the maximal ideal $\fom$.
Parameter ideals have the advantage that one may consider
the tight closure question whether
$f \in (f_1 \komdots f_d)^*$ holds as a problem about the $\rm \check{C}$ech cohomology class
$\frac{f}{f_1 \cdots f_d} \in H^d_\fom(R)$.

The general question remained open even for homogeneous $R_+$-primary
ideals in a two-dimensional normal
standard-graded domain.
After a long period without substantial progress
it was shown in \cite{brennertightelliptic} that the question has a positive
answer for homogeneous primary ideals in affine cones over an elliptic curve,
such as $R=K[X,Y,Z]/(X^3+Y^3+Z^3)$, itself a prominent example in tight closure theory.
This result rests on a geometric interpretation of tight closure theory
in terms of vector bundles and uses the classification of
vector bundles on elliptic curves due to M. Atiyah.

In this paper we show that $I^* =I^+$ (in fact we even show that
$I^*=I^{\gr}$, the graded plus closure) holds for a homogeneous
ideal in a two-dimensional $\NN$-graded domain $R$ of finite type
over the algebraic closure of a finite field (Theorem
\ref{tightplusclosure2}). This last finiteness condition is not
necessary in the elliptic case, but it is crucial for our proof for
higher genus. The point is that due to a result of H. Lange and U.
Stuhler \cite{langestuhler} a strongly semistable bundle (see below
for the definition) of degree $0$ on a smooth projective curve $Y$
defined over a finite field can be trivialized by pulling it back
along a finite mapping $Y' \ra Y$.

The other main ingredient of our proof is \---beside the geometric
interpretation of tight closure\--- a recent theorem of A. Langer.
He shows for a locally free sheaf on a smooth projective variety
that the Harder-Narasimhan filtration of some Frobenius pull-back
has strongly semistable quotients \cite[Theorem
2.7]{langersemistable}. This allows us to do induction on this
strong Harder-Narasimhan filtration and to obtain for the tight
closure and the plus closure (or rather for their geometric
counterparts) the same numerical criterion.

The paper is organized as follows. In Section \ref{pre} we describe
briefly the geometric interpretation and the two mentioned results
which we will use in the following.

Sections \ref{degreecriteria} and \ref{slopetrivial}
deal with the geometric setting
which arises from tight closure:
a locally free sheaf $\shS$ on a smooth projective curve $Y$
and a geometric torsor $T \ra Y$ corresponding to a cohomology class $c \in H^1(Y,\shS)$.
This setting is of independent interest and no knowledge of tight closure
is required.

Section \ref{degreecriteria} establishes a numerical criterion for
the affineness of such a torsor on a curve in every characteristic
(Theorem \ref{numcritaffine}). Section \ref{slopetrivial}
establishes the same numerical criterion for the existence of a
projective curve inside the torsor under the condition that the
curve is defined over a finite field (Theorem
\ref{numcritprojcurve}).

In Section \ref{plustight} we derive the consequences from our
results to tight closure. The equality $I^*=I^{\gr}$ for
$R_+$-primary homogeneous ideals in the normal standard-graded case
follows immediately (Theorem \ref{tightplusclosure}) and in Theorem
\ref{tightplusclosure2} we remove the conditions primary, normal and
standard-graded.

\section{Preliminaries}
\label{pre}

We describe briefly the geometric setting of our approach to tight
closure and plus closure (see \cite{brennertightproj} for details).
Let $R$ denote a normal standard-graded domain over an algebraically
closed field $K$ (of any characteristic) and let $f_1 \komdots f_n$
denote homogeneous generators of an $R_+$-primary ideal of degrees
$d_1 \komdots d_n$. These data give rise to the short exact sequence
of locally free sheaves on $Y= \Proj R$,
$$ 0 \lra \Syz(f_1 \komdots f_n)(m) \lra
 \bigoplus_{i=1}^n  \O(m-d_i) \stacklra {f_i}  \O(m) \lra 0 \, .$$
Another homogeneous element $f \in R$ of degree $m$ yields via the connecting
homomorphism
the cohomology class $c=\delta(f) \in H^1(Y, \Syz(f_1 \komdots f_n)(m))$.
Such a class $c \in H^1(Y,\shS)$ of a locally free sheaf $\shS$
yields a geometric $\shS$-torsor $T \ra Y$,
that is an affine-linear bundle on which $\shS$ acts.

With this setting, the containment $f \in (f_1 \komdots f_n)^*$ is
in positive characteristic equivalent to the property that the
cohomological dimension of the torsor $T$ equals $d=\dim Y =\dim
R-1$ \cite[Proposition 3.9]{brennertightproj}. If $R$ has dimension
two, this just means that $T$ is not an affine scheme.

The graded plus closure of an ideal, denoted $I^{\gr}$, consists of
the elements such that there exists a finite homogeneous extension
$R \subseteq S$ (of graded domains) such that $f \in IS$. In terms
of our geometric setting, the containment $f\in (f_1 \komdots
f_n)^{\gr} $ is equivalent to the geometric property that the torsor
$T$ contains a closed projective subvariety of dimension $d$
\cite[Lemma 3.10]{brennertightproj}. This is equivalent to the
property that there exists a finite dominant mapping $\varphi: Y'
\ra Y$ of a projective normal variety $Y'$ such that
$\varphi^*(c)=0$ in $H^1(Y', \varphi^*(\shS))$.

We restrict now to the case of a two-dimensional normal domain $R$
so that $Y$ is a smooth projective curve. A crucial notion to study
the affineness of a torsor $T$ given by a cohomology class $c \in
H^1(Y, \shS)$ is the semistability of $\shS$. Recall that a locally
free sheaf $\shS$ is called semistable if $\mu(\shT) \leq \mu(\shS)$
holds for every subbundle $\shT \subseteq \shS$, where $\mu(\shT) =
\deg (\shT)/ \rk (\shT)$ denotes the slope of a bundle (see
\cite{huybrechtslehn} for background of this notion).

In positive characteristic, the pull-back of a semistable bundle
under the absolute Frobenius $F: Y \ra Y$ is in general not
semistable. However, if it stays semistable for every Frobenius
power, then the bundle is called strongly semistable, a notion
introduced by Miyaoka in \cite{miyaokachern}. Moreover, a strongly
semistable bundle stays semistable under every finite mapping
$\varphi:Y' \ra Y$ \cite[Proposition 5.1]{miyaokachern}.

A useful characterization for a strongly semistable bundle on a
curve defined over a finite field is given in \cite{langestuhler}.
In this paper, H. Lange and U. Stuhler study $\GL$-representations
of the algebraic fundamental group in positive characteristic and
stability properties of the corresponding vector bundles. One of
their results says that a vector bundle $\shS$ becomes trivial on an
\'{e}tale covering (and corresponds then to a continuous
re\-presentation of the algebraic fundamental group $\pi_1(Y) \ra
\GL(r,K)$) if and only if some Frobenius pull-back of $\shS$ is
isomorphic to $\shS$ (which is only possible for $\deg(\shS)=0$).
For our purpose the following result is more important.

\begin{theorem}
\label{finitetrivialstrongly}
Let $K$ denote the algebraic closure of a finite field of characteristic $p$
and let $Y$ denote a smooth projective curve over $K$.
Let $\shS$ denote a locally free sheaf on $Y$ of degree $0$.
Then the following are equivalent.

\numiii

\begin{enumerate}

\item
$\shS$ trivializes under a finite mapping $Y' \ra Y$.

\item
$\shS$ is strongly semistable.

\item
There exist numbers $e'> e$ such that $F^{e' *} (\shS) \cong F^{e*}(\shS)$.

\item
There exists $Y' \stackrel {\varphi}{\lra} Y \stackrel{F^{e}}{\lra}
Y$, where $\varphi$ is \'{e}tale and where the pull-back
$(F^{e}\circ \varphi)^*(\shS)$ is trivial.
\end{enumerate}
\end{theorem}
\proof This is proven in \cite[Satz 1.9, Korollar 1.6 and Satz
1.4]{langestuhler}. \qed

\begin{remark}
The pull-back of a syzygy bundle $\Syz(f_1 \komdots f_n)(m)$ under
the $e$-th Frobenius morphism is given by $\Syz(f_1^q \komdots
f_n^q)(qm)$, where $q=p^{e}$. So the condition (iii) in Theorem
\ref{finitetrivialstrongly} translates to $\Syz(f_1^q \komdots
f_n^q)(qm) \cong \Syz(f_1^{q'} \komdots f_n^{q'})(q'm)$. If moreover
$\Syz(f_1^q \komdots f_n^q)(qm) \cong \Syz(f_1 \komdots f_n)(m)$ for
some $q=p^{e}$, then there exists even an \'{e}tale covering $Y'$ of
$Y$ such that this syzygy bundle becomes trivial.
\end{remark}

For every locally free sheaf $\shS$ on the smooth projective curve
$Y$ there exists the so-called Harder-Narasimhan filtration $\shS_1
\subsetdots \shS_t =\shS$. This filtration is unique and has the
property that the quotients $\shS_{i+1}/ \shS_{i}$ are semistable
and $\mu (\shS_{i}/ \shS_{i-1}) > \mu(  \shS_{i+1}/ \shS_{i})$ holds
true. The number $\mu(\shS_1) = \mu_{\rm max}(\shS)$ is called the
maximal slope of $\shS$ and $\mu(\shS/ \shS_{t-1}) = \mu_{\rm
min}(\shS)$ the minimal slope of $\shS$.

The pull-back of the Harder-Narasimhan filtration under the
Frobenius yields in general not the  Harder-Narasimhan filtration of
$F^*(\shS)$. However, a recent result of A. Langer \cite[Theorem
2.7]{langersemistable} shows that there exists a Frobenius power
$F^{e}$ such that the quotients in the Harder-Narasimhan filtration
of the pull-back $F^{e*}(\shS)$ are all strongly semistable. We call
such a filtration the strong Harder-Narasimhan filtration.

This strong Harder-Narasimhan filtration is useful for us in several respects.
It allows us to
reduce some questions to the strongly semistable case, where
a numerical characterization for tight closure was given in \cite[Theorem 8.4]{brennerslope}.

It also implies that the number
$$\bar{\mu}_{\min}(\shS) = \lim_{e \in \NN}  \frac{\mu_{\rm min} (F^{e *}(\shS))}{p^{e}}
=\min \{ \frac{\mu_{\min} (\varphi^*(\shS))}{\deg (\varphi)}:\,\, \varphi : Y' \ra Y \mbox{ finite} \}
$$
is a rational number and that it is obtained for some $e$. This
improves also the slope criterion for ample vector bundles: $\shS$
is ample if and only if $\bar{\mu}_{\min}(\shS) >0$ \cite[Theorem
2.3]{brennerslope}, and this is now equivalent to $ \mu_{\min}(
F^{e*}\shS) >0$ for all $e \in \NN$.

Furthermore, if $\shS_1 \subsetdots \shS_t =\shS$ is the strong
Harder-Narasimhan filtration, then we may look for the maximal $i$
such that $\mu( \shS_i/ \shS_{i-1}) \geq 0$ (if $\mu( \shS_j/
\shS_{j-1}) < 0$ for all $j=1 \komdots t$, then set $i=0$ and
$S_0=S_{-1}=0$ and if $\mu( \shS_j/ \shS_{j-1}) \geq 0$ for all $j=1
\komdots t$, then set $i=t+1$ and $\shS_{t+1}=\shS_t=\shS$). It will
be crucial to consider the short exact sequence $0 \ra \shS_i \ra
\shS \ra \shS/\shS_i = \shQ \ra 0$. Of course, $\shS_1 \subsetdots
\shS_i$ is the strong Harder-Narasimhan filtration of $\shS_i$ and $
\shS_{i+1}/\shS_i \subset \shS_{i+2}/\shS_i \subsetdots \shS/\shS_i
=\shQ$ is the strong Harder-Narasimhan filtration of $\shQ$.
Therefore $\bar{\mu}_{\rm min} (\shS_i) \geq 0$ and $\bar{\mu}_{\rm
max}(\shQ) <0$. Hence $\bar{\mu}_{\min}(\shQ^\dual) >0$ and the dual
$\shQ^\dual$ is an ample vector bundle.

\section{Slope criteria for the affineness of torsors}
\label{degreecriteria}

Let $Y$ be a smooth projective curve over an algebraically closed field $K$ of
arbitrary characteristic.
Let $\shS$ denote a locally free sheaf on $Y$.
A cohomo\-lo\-gy class $c \in H^1(Y,\shS)$ corresponds
to a geometric torsor $T \ra Y$.
This is an affine-linear bundle on which $\shS$ acts by translations.
A natural realization of $T$ is obtained as follows:
since $H^1(Y,\shS) \cong \Ext(\O_Y, \shS)$,
the class yields an extension
$0 \ra \shS \ra \shS' \ra \O_Y \ra 0$
which gives a projective embedding $\PP(\shS^\dual) \hookrightarrow \PP(\shS^{'\dual})$.
Then $T=\PP(\shS^{'\dual}) - \PP(\shS^{\dual})$.

We are interested in the global properties of such a torsor $T \ra
Y$ corresponding to $c \in H^1(Y,\shS)$. In this section we give a
numerical criterion in terms of the strong Harder-Narasimhan
filtration in order to decide whether such a torsor is an affine
scheme or not.

\begin{proposition}
\label{ampleaffine}
Let $Y$ denote a smooth projective curve over an algebraically closed field $K$.
Let $\shS$ denote a locally free sheaf and let $c\in H^1(Y,\shS)$
be a cohomology class with corresponding torsor $T \ra Y$.
If $\shS$ is strongly semistable, then the following are equivalent.

\numiii

\begin{enumerate}
\item
The torsor $T$ is an affine scheme.

\item
$\bar{\mu}_{\max}(\shS) < 0$ and $c \neq 0$ {\rm(}in positive
characteristic $F^{e*}(c) \neq 0$ for all Frobenius powers
$F^{e}${\rm)}.

\item
The dual extension $\shS^{'\dual}$ given by $c$ is an ample bundle.
\end{enumerate}
The implications {\rm (ii)} $\Leftrightarrow$ {\rm (iii)} $\Rightarrow$ {\rm (i)}
hold for any locally free sheaf $\shS$.
\end{proposition}
\proof (i) $\Rightarrow $ (ii). Since $T$ is affine, we have $c \neq
0$, for otherwise $T$ would be trivial and contain projective
curves. Since the pull-back $T'=T \times_Y Y'$ for $Y' \ra Y$ finite
is again affine, we get also $F^{e*}(c) \neq 0$ in positive
characteristic. The affineness of $T$ implies in general
$\bar{\mu}_{\min}(\shS)<0$ due to \cite[Theorem 4.4]{brennerslope}.
But if $\shS$ is strongly semistable, then $\bar{\mu}_{\min}(\shS)
= \bar{\mu}_{\max}(\shS)<0$.

(ii) $\Rightarrow $ (iii). The slope condition yields
$\bar{\mu}_{\min}(\shS^\dual) > 0$ for the dual bundle $\shS^\dual$.
Hence $\shS^\dual$ is an ample vector bundle (see the discussion at
the end of the previous section). Since $\shS^{'\dual}$ is a
non-trivial extension, it follows from \cite[Proposition
2.2]{giesekerample} that all quotients of $\shS^{'\dual}$ have
positive degree. Since this is true for every Frobenius pull-back,
we have again $\bar{\mu}_{\min}(\shS^{'\dual}) > 0$ and so
$\shS^{'\dual}$ is also ample.

(iii) $\Rightarrow $ (i).
The ampleness of $\shS^{'\dual}$ is by definition the ampleness of the divisor
$\PP(\shS^\dual) \subset \PP(\shS^{'\dual})$.
So the complement of $\PP(\shS^\dual)$ is an affine scheme.

Finally we prove (iii) $\Rightarrow $ (ii) without the condition
that $\shS$ is strongly semi\-stable. The surjection $\shS'^\dual
\ra \shS^\dual \ra 0$ shows that $\shS^\dual$ is also ample; hence
$\bar{\mu}_{\min}(\shS^\dual) > 0$ and therefore
$\bar{\mu}_{\max}(\shS)<0$. If $F^{e*}(c)=0$, then
$F^{e*}(\shS'^\dual) \cong  \shS^\dual \oplus \O_Y$ would be not
ample. \qed

\begin{remark}
It may indeed happen that a class $0 \neq c \in H^1(Y,\shS)$ for
$\shS$ of negative degree (even for invertible $\shS$) becomes zero
under a Frobenius power. But this is an exception. In the case of a
relative curve over $\Spec \ZZ$, this may happen only for finitely
many prime numbers. For the condition $c \neq 0$ implies ampleness
over the generic curve of characteristic zero, and this holds then
almost everywhere, since ampleness is an open property.
\end{remark}

We look now at the affineness of the torsor $T$ for an arbitrary
locally free sheaf $\shS$. The crucial point is to look at the
strong Harder-Narasimhan filtration of $\shS$ (in characteristic
zero this is just the usual Harder-Narasimhan filtration, and
replace the Frobenius by the identity). Note that the property of
being affine is preserved under a finite mapping. So in positive
characteristic we often apply an absolute Frobenius morphism $F:Y\ra
Y$ to pull back the whole situation.

\begin{theorem}
\label{numcritaffine} Let $\shS$ denote a locally free sheaf on a
smooth projective curve $Y$ over an algebraically closed field $K$,
let $c\in H^1(Y,\shS)$ with corresponding torsor $T\ra Y$. Let
$\shS_1 \subsetdots \shS_t $ be the strong Harder-Narasimhan
filtration of $F^{e*} (\shS)$ on $Y$. Choose $i$ such that
$\shS_i/\shS_{i-1}$ has degree $\geq 0$ and that $\shS_{i+1}/
\shS_i$ has degree $<0$ {\rm(}we set $i=-1$ and $\shS_{-1}=0$ and
$i=t+1$ and $\shS_{t+1}=F^{e*} (\shS)$ in the extremal cases{\rm)}.
Set $0 \ra \shS_i \ra F^{e*}(\shS) \ra F^{e*}(\shS)/\shS_i =\shQ \ra
0$. Then $T$ is affine if and only if the image of $F^{e*}(c)$ in
$H^1(Y, \shQ)$ is non-zero {\rm(}in positive characteristic non-zero
for all Frobenius powers{\rm)}.
\end{theorem}
\proof Let $c' \in H^1(Y, \shQ)$ denote the image of $F^{e*}(c)$ and
suppose first that $c'\neq 0$. By construction, $\shQ$ is a locally
free sheaf $\neq 0$ with $\bar{\mu}_{\rm max}(\shQ)<0$. So its dual
sheaf $\shQ^\dual$ is an ample vector bundle and the same holds true
for the extension $\shQ^{'\dual}$ defined by $c'$ as in the proof of
Proposition \ref{ampleaffine}. Therefore the torsor $T'$
corresponding to $c'$ is affine by Proposition \ref{ampleaffine}.
There exists an affine morphism $F^{e*}(T) \ra T'$ induced by
$F^{e*}(\shS) \ra \shQ$ \cite[Lemma 3.1]{brennertightelliptic},
hence $F^{e*}(T)$ and then also $T$ itself is affine.

On the other hand, suppose that some Frobenius pull-back of the
image $c'$ is $0$. Then we may assume $c'=0$ and hence $F^{e*}(c)$
stems from a cohomology class $\tilde{c} \in H^1(Y, \shS_i)$, where
$\bar{\mu}_{\rm min} (\shS_i) \geq 0$ (including the case
$\shS_i=\shS_{-1}=0$). It is enough to show that the affine-linear
bundle $\tilde{T}$ corresponding to $\tilde{c}$ is not affine. But
this was proven in \cite[Theorem 4.4]{brennerslope}. \qed

\section{Slope criteria for the trivialization of a cohomology class}
\label{slopetrivial}

In this section we suppose that $K$ is the algebraic closure of a
finite field of positive characteristic $p$, $K=\overline{\FF}_p$.
We shall now prove for a smooth projective curve $Y$ defined over
$K$ that the criterion of the last section for the affineness of a
torsor holds also for the property that there does not exist any
projective curve inside the torsor. This property itself is
equivalent to the fact that the cohomolgy class does not get trivial
under any finite mapping $\varphi:Y' \ra Y$. We first do the case of
a strongly semistable bundle.

\begin{theorem}
\label{strongsemistabletrivial}
Let $Y$ denote a smooth projective curve over the algebraic closure of a finite field.
Let $\shS$ denote a strongly semistable locally free sheaf of degree $\geq 0$ and let
$c \in H^1(Y, \shS)$ denote a cohomology class.
Then $c$ trivializes under a finite mapping $Y' \ra Y$.
\end{theorem}
\proof Suppose first that the degree of $\shS$ is positive. Then
$\shS$ is an ample bundle and so it is also cohomologically p-ample
in the sense of \cite{kleimanamplesurface} (see
\cite{migliorinicohomological} or \cite{arapuraamplitude}), that is,
for every coherent sheaf $\shF$ we have $H^{i}(Y, \shF \otimes
F^{e*}(\shS)) =0$ for $i \geq 1$ and $e$ large enough. This implies
in particular that some Frobenius power of the class $c$ is trivial.

So suppose that $\deg (\shS)=0$. Since $K$ is supposed to be the
closure of a finite field, every object we encounter is in fact
defined already over a finite field. Since $\shS$ is strongly
semistable, there exists by the result of Lange and Stuhler (see
Theorem \ref{finitetrivialstrongly} above) a finite mapping
$\varphi: Y' \ra Y$ (which is an \'{e}tale mapping followed by a
Frobenius power) such that the pull-back is trivial, hence
$\varphi^*(\shS) = \O_{Y'}^r$. Thus we may assume that $\shS$ is in
fact a trivial bundle. We may deal with the components of $c \in
H^1(Y, \O_Y^r)$ separately, so we may even assume that $\shS= \O_Y$.
But for the structure sheaf this was proven in \cite[Proposition
8.1]{brennertightproj} and also follows from \cite[Proposition
3.3]{smithgraded}. \qed

\begin{theorem}
\label{numcritprojcurve} Let $Y$ denote a smooth projective curve
over the algebraic closure of a finite field. Let $\shS$ denote a
locally free sheaf on $Y$ and let $c\in H^1(Y,\shS)$ with
corresponding torsor $T \ra Y$. Let $\shS_1 \subsetdots \shS_t $ be
the strong Harder-Narasimhan filtration of $F^{e*}(\shS)$ on $Y$.
Choose $i$ such that $\shS_i/\shS_{i-1}$ has degree $\geq 0$ and
that $\shS_{i+1}/ \shS_i$ has degree $<0$. Let $0 \ra \shS_i \ra
F^{e*}(\shS) \ra F^{e*}(\shS)/\shS_i =\shQ \ra 0$. Then the
following are equivalent:

\numiii

\begin{enumerate}

\item
The $\shS$-torsor $T$ contains a projective curve.

\item
There exists a smooth projective curve $Y'$ and a finite mapping
$\varphi: Y' \ra Y$ such that $\varphi^*(c)=0$.

\item
Some Frobenius power of the image of $F^{e*}(c)$ in $H^1(Y, \shQ)$ is zero.

\end{enumerate}

\end{theorem}

\proof
The equivalence of (i) and (ii) is clear, since a torsor is trivial
if and only if it has a section.
Suppose that (iii) holds.
We may assume that the image of $F^{e*}(c)$ is $0$ inside $H^1(Y, \shQ)$.
Then $F^{e*}(c)$ stems from a class $c_i \in H^1(Y, \shS_i)$.
The filtration $0 \subset S_1 \subsetdots S_{i-1} \subset S_i$
is such that all quotients are strongly semistable
with $\deg(\shS_j/ \shS_{j-1}) \geq 0$ for $j=1 \komdots i$.
We look at the short exact sequence
$0 \ra \shS_{i-1} \ra \shS_i \ra \shS_i/\shS_{i-1} \ra 0$.
Since the sheaf on the right is strongly semistable of non-negative degree,
it follows from Theorem \ref{strongsemistabletrivial}
that the image $c_i'$ of $c_i$ in $H^1(Y, \shS_i/\shS_{i-1})$
vanishes under a finite mapping $\varphi: Y' \ra Y$.
Then $\varphi^*(c_i)$ stems from a class
$c_{i-1} \in H^1(Y', \varphi^*(\shS_{i-1}))$.
Going on like this inductively we find a finite mapping such that the pull-back of $c$
is zero.

Condition (i) implies that $T$ is not an affine scheme.
So (i) $\Rightarrow $ (iii) follows from Theorem \ref{numcritaffine}.
\qed

\medskip
Bringing the results of this and the last section together we deduce the following
corollary, which gives a geometric criterion
for the cohomological property of being affine.

\begin{theorem}
\label{affineprojectivecurve}
Let $Y$ denote a smooth projective curve over the algebraic closure of a finite field.
Let $\shS$ denote a locally free sheaf on $Y$
and let $c\in H^1(Y,\shS)$ with corresponding geometric torsor $T \ra Y$.
Then $T$ is an affine scheme if and only if it does not contain any projective curve.
\end{theorem}

\proof
This follows from Theorem \ref{numcritaffine} and Theorem \ref{numcritprojcurve},
since for both properties the same numerical criterion holds.
\qed

\section{Plus closure and tight closure}
\label{plustight}

We come now back to tight closure. The results in the previous sections
give at once the identity $I^*=I^+$ for homogeneous
primary ideals in a normal standard-graded two-dimensional domain.

\begin{theorem}
\label{tightplusclosure} Suppose that $K$ is the algebraic closure
of a finite field. Let $R$ denote a standard-graded two-dimensional
normal domain over $K$. Then for every homogeneous $R_+$-primary
ideal $I$ we have the identities $ I^* =I^{\gr } =I^+$.
\end{theorem}
\proof
The inclusions $I^{\gr} \subseteq I^+ \subseteq I^*$ are clear,
so it is enough to show that $I^* \subseteq I^{\gr }$.
Due to \cite[Theorem 4.2]{hochsterhunekesplitting}
we only have to consider homogeneous elements.
Let $I=(f_1 \komdots f_n)$ be given by homogeneous ideal generators.
These give rise to the locally free syzygy bundle
$\Syz(f_1 \komdots f_n)(m)$ on the smooth projective curve $Y =\Proj\, R$.
A homogeneous element $f$ of degree $m$ induces a cohomology class
$c=\delta (f) \in H^1(Y,\Syz(f_1 \komdots f_n)(m))$ with corresponding
torsor $T \ra Y$.
The containment of $f$ in the graded plus closure,
$f \in I^{\gr }$, is equivalent to the existence of a projective curve inside $T$.
The containment in the tight closure,
$f \in I^*$, is equivalent to the non-affineness of $T$.
So the result follows from Theorem \ref{affineprojectivecurve}.
\qed

\begin{theorem}
\label{tightplusclosure2}
Suppose that $K$ is the algebraic closure of a finite field.
Let $R$ denote an $\NN$-graded two-dimensional domain of finte type over $K$.
Then for every homogeneous ideal $I$ we have the identities
$$ I^* =I^{\gr } =I^+  \, .$$
\end{theorem}
\proof
We first reduce the identity to the primary case.
Let $I=(f_1 \comdots f_n)$ denote a homogeneous ideal
and suppose $f \in I^*$, $f$ homogeneous of degree $m$.
Then for every $k \in \NN$
we have $f \in (I +R_{\geq k})^*$. Since these are $R_+$-primary ideals, we have
also $f \in (I +R_{\geq k})^{\gr}$.
This means that for every $k$ we have a finite, homogeneous (degree-preserving)
extension $R \ra S$, where $S$ is another graded domain,
with $f \in (I+R_{\geq k})S$.
This means that we have an equation
$f = \sum_{i=1}^ns_if_i + \sum_{j} t_jg_j$, where
$s_i, t_j \in S$ and $g_j \in R_{\geq k}$.
We may assume that everything is homogeneous,
hence for $k >m$ we get $t_j=0$ and therefore $ f \in IS$.

Now suppose that $R$ is a two-dimensional $\NN$-graded domain $R$ of finite
type over the algebraic closure $K$ of a finite field.
Write $R=K[T_1 \komdots T_k]/\foa$ with $\deg (T_i)=e_i$,
and look at the homogeneous ring homomorphism
$$K[T_1 \komdots T_k] \ra K[U_1 \komdots U_k],\, \, T_i \mapsto U^{e_i} $$
and the induced mapping $R \ra K[U_1 \komdots U_k]/ \foa K[U_1 \komdots U_k]=:R'$.
$R'$ is now a standard-graded two-dimensional $K$-algebra finite
over $R$. We may mod out a homogeneous minimal prime ideal of $R'$ and normalize
to get a finite mapping $R \ra S$,
where $S$ is now a normal two-dimensional standard-graded domain.
From $f \in I^*$ we get $f \in (IS)^*$ and hence due
to Theorem \ref{tightplusclosure} also $f \in (IS)^{\gr }$.
But then also $f \in I^{\gr}$.
\qed

\medskip
Of course it is natural to ask whether Theorems \ref{tightplusclosure2}
and Theorem \ref{affineprojectivecurve}
hold without the assumption that $K$ is the algebraic closure of a finite field.
A special case of this question is the following problem.

\begin{problem}
Let $\shL$ denote an invertible sheaf of degree $0$ on a smooth projective curve
$Y$ over an algebraically closed field $K$ of positive characteristic.
Let $c \in H^1(Y, \shL)$. Does there exist a finite
mapping $\varphi: Y' \ra Y$ such that $\varphi^*(c)=0$?

This is true for the structure sheaf and for every $\shL$ which is a torsion point in the
Picard group. Therefore it is true in general for the algebraic closure
of a finite field.
\end{problem}

\bibliographystyle{plain}

\bibliography{bibliothek}

\end{document}